\documentclass[]{amsart}
\author{K. Denecke, J. Koppitz, Sl. Shtrakov}
\title{\bf The Depth of a Hypersubstitution}
\date{}

\begin{document}

\maketitle
\sloppy
\newfont{\zb}{cmssi12 scaled 1000}
\newcommand{\natur}{\mbox{\it{I\hspace{-0.2em}N}}}

\newcounter{nummer}
\newtheorem{satz}{Theorem}[section]
\newtheorem{df}[satz]{Definition}
\newtheorem{bsp}[satz]{Example}
\newtheorem{lemma}[satz]{Lemma}
\newtheorem{leer}[satz]{}
\newtheorem{prop}[satz]{Proposition}
\newtheorem{cor}[satz]{Corollary}

\vspace{1cm}

\begin{abstract}
For given depths of the terms $s, t_1, \cdots, t_n$  a formula will be proved
to calculate the depth of the composed term $s(t_1, \cdots, t_n)$ and
if $\sigma$ is a hypersubstitution and $t$ is a term we derive a formula for
the depth of $\hat \sigma[t]$.

\end{abstract}

\section{Introduction}
At first we remember of the following definition of terms.
Let $X = \{x_1, \cdots, x_n, \cdots\}$ be any countably infinite (standard)
alphabet of variables and let $X_n = \{x_1, \cdots, x_n\}$ be an $n$-element
alphabet. Let $(f_i)_{i \in I}$ be an indexed set which is disjoint from
$X$. Each $f_i$ is called an $n_i$-ary operation symbol where $n_i \geq 1$ is a
natural number.  Let $\tau$ be a function which assigns to every $f_i$ the
number $n_i$ as its arity. The function $\tau$ or the sequence of values
of $\tau$, written as $(n_i)_{i \in I}$, is called a type. An $n$-ary term of
type $\tau$ is defined inductively as follows:
\begin{enumerate}
\item[\mbox{(i)}] The variables $x_1, \cdots, x_n$ are $n$-ary terms.
\item[\mbox{(ii)}] If $t_1, \cdots, t_m$ are $n$-ary terms and if $f_i$
is an $n_i$-ary operation symbol then $f_i(t_1, \cdots, t_{n_i})$ is an $n$-ary
term.
\item[\mbox{(iii)}] Let $W_{\tau}(X_n)$ be the smallest set which contains
$x_1, \cdots, x_n$ and is closed under finite application of (ii). Every
$t \in W_{\tau}(X_n)$ is called an $n$-ary term of type $\tau$.
\end{enumerate}
We remark that by this definition every $n$-ary term is also $(n+1)$-ary.
The set $W_{\tau}(X):= \bigcup \limits_{n=1}^{\infty}W_{\tau}(X_n)$ is the
set of all terms of type $\tau$.

Usually one has a third set $\overline A$, called set of constants, with
$\overline A \cap X = \emptyset$ and $\{f_i|i \in I\}
\cap (\overline A \cup X) = \emptyset.$
Then polynomials of type $\tau$ over $\overline A$ are defined in a
similar way adding a condition which says that constants are polynomials.
Constants can also be defined by nullary operation symbols assuming that
the indexed set $(f_i)_{i \in I}$ of operation symbols includes also nullary
operation symbols $(n_i = 0)$.
We mention also that in the case of finite sets of operation symbols instead of
polynomials one speaks of trees. Trees can be regarded as connected
graphs without cycles. Trees have many applications in Computer Science, in
Linguistic and in other fields.
Ordered binary decision diagrams (OBDD's) are trees in the language of
Boolean algebras. Trees can be used to visualize the structure of computer
programmes. For all these applications it is important to measure the
complexity of a tree or of a term. The concept of the depth of a term
(or of the height of a tree) is a well-known complexity measure.
The method of algebraic induction which is very often used, is based on the
depth of terms or polynomials.

\begin{df}
Let $t \in W_{\tau}(X)$ be a term.
\begin{enumerate}
\item[\mbox{(i)}] If $t= x\in X$ then $Depth(t):= 0$.
\item[\mbox{(ii)}] If $t = f_i(t_1, \cdots, t_{n_i})$ then\\ $Depth(t):=
max\{Depth(t_1), \cdots, Depth(t_{n_i})\} + 1$.
\end{enumerate}
\end{df}
The depth of a polynomial is defined in a similar way where
$Depth(\overline a):= 0$ if $\overline a$ is a constant from $\overline A$.\\

In the case of a tree usually one speaks of the height of the tree.\\

Our goal is to describe the behavior of the depth under
some mappings defined on sets of terms. We select two mappings which play
an important role in Universal Algebra but also in Computer Science.\\

The first mapping is called composition of terms and is defined in the
following inductive way:\\
Let $s \in W_{\tau}(X_n)$ and let $t_1, \cdots, t_n \in W_{\tau}(X_m)$.
Then we define
\[S_m^n: W_{\tau}(X_n) \times W_{\tau}(X_m)^n \to W_{\tau}(X_m)\]
by the following steps:
\begin{enumerate}
\item[\mbox{(i)}] If $s=x_i, 1 \leq i \leq n$
then $S_m^n(s, t_1, \cdots, t_n):= t_i$.
\item[\mbox{(ii)}] If $s = f(s_1, \cdots, s_r)$ and $s_1, \cdots, s_r \in
W_{\tau}(X_n)$ then
\item[] $S_m^n(s,t_1, \cdots, t_n):= f(S_m^n(s_1, t_1, \cdots, t_n), \cdots,
S_m^n(s_r, t_1, \cdots, t_n))$
\item[] $= f(s_1(t_1, \cdots, t_n),
\cdots, s_r(t_1, \cdots, t_n))$.
\end{enumerate}

We remark that the heterogeneous (multibased) algebra
$$(W_{\tau}(X_n)_{n \in \natur^+}, (S_m^n)_{m,n \in \natur^+},
(x_i)_{i \leq n \in \natur^+})$$
with $\natur^+ :=\natur \setminus \{0\}$ with variables as nullary operations
is called full term clone of type $\tau$. (There are only technical reasons not
to consider nullary terms, but one can define a similar structure
for polynomials as well.)
For proofs by induction it is very advantageous to use the
operations $S_m^n$ to describe the superposition of terms.

The second kind of mappings is called hypersubstitution of type $\tau$, for
short, hypersubstitution.
Hypersubstitutions play an important role in the theory of hyperidentities
and solid varieties which is a very fast-developing modern algebraic
theory with applications in Theoretical Computer Science (hyper-tree automata)
and in Logic (fragment of second order logic).

Hypersubstitutions $\sigma$ of type $\tau$ are defined by
$\sigma(f_i) = t$ where $t$ is an $n_i$-ary term for all $n_i-$ary
operation symbols $f_i$.
Those mappings can be extended to mappings $\hat \sigma$ defined on sets of
terms by the following steps:\\
\begin{enumerate}
\item[\mbox{(i)}] $\hat \sigma[x]:= x$ if $x \in X$ is a variable.
\item[\mbox{(ii)}] $\hat \sigma[f_i(t_1, \cdots, t_{n_i})]:=
S_n^{n_i}(\sigma(f_i),\hat \sigma[t_1], \cdots, \hat \sigma[t_{n_i}])$ for composed
terms $f_i(t_1, \cdots, t_{n_i})$.
\end{enumerate}
This definition shows that the mappings $\hat \sigma$ are endomorphisms of the
full term clone (as heterogeneous algebra). In the case of trees
the mappings $\hat \sigma$ are special types of so-called alphabetic tree
homomorphisms (\cite{Gec, Com}) and preserve the
recognizability of a forest (set of trees) by a tree automaton.\\
Having a look on trees the formula which we will derive for the depth
of composed trees becomes quite clear and simple.
Nevertheless we will give full proofs, mostly by induction since terms
can be countably infinite. For the depth of an
arbitrary hypersubstitution we obtained a more complicated formula and its
proof is far from beeing trivial.

\section{Full Terms}

Our first aim is to calculate $Depth(S_m^n(s,t_1, \cdots, t_n))$ if\\ $Depth(s),
Depth(t_1), \cdots, Depth(t_n)$ are known.
The following example shows that $Depth(S_m^n(s,t_1, \cdots, t_n))$ depends
not only on the depths of the inputs but also on the special structure
of the term $S_m^n(s,t_1, \cdots, t_n)$.

Consider the type $\tau =(2)$ with $f$ as binary operation symbol and the terms
$$t_1 = f(x_1, f(x_1,x_2)), t_2 = f(x_2,x_1)\ \ \mbox{ and }\ \  s_1 = f(f(x_2,x_2),x_1).$$
Then we have  $$Depth(t_1) = 2, Depth(t_2) = 1, Depth(s_1) = 2\ \ \mbox{ and }$$
$$S_2^2(s_1,t_1,t_2) = s_1(t_1,t_2) = f(f(f(x_2,x_1), f(x_2,x_1)),
f(x_1, f(x_1,x_2))),\ \ \mbox{ and }$$ $$ Depth(S_2^2(s_1,t_1,t_2)) = 3.$$
Now instead of $s_1$ we take the term $s_2 = f(f(x_1,x_1),x_2)$ with
$Depth(s_2) = 2$ and obtain $$S_2^2(s_2,t_1,t_2) =  f(f(f(x_1,f(x_1,x_2)),
f(x_1,f(x_1,x_2))), f(x_2,x_1))$$ with  $Depth(S_2^2(s_2,t_1,t_2)) = 4$.
That means, $Depth(S_m^n(s, t_1, \cdots, t_n))$ depends on the particular
structure of the terms. But this is not always the case. If, for instance,
$s = x_i, 1 \leq i \leq n$ is a variable then $Depth(s) = 0$ and
$Depth(S_m^n(s,t_1, \cdots, t_n)) = Depth(t_i)$.

Now we consider the following kind of terms, called full terms:
\begin{df}\label{d2.1}
\begin{enumerate}
\item[\mbox{(i)}] If $f_i$ is an $n_i$-ary operation symbol and if $s:
\{1, \cdots, n_i\} \to \{1, \cdots, n_i\}$ is a permutation then
$f_i(x_{s(1)}, \cdots, x_{s(n_i)})$ is a full term.
\item[\mbox{(ii)}] If $f_j$ is an $n_j$-ary operation symbol and if
$t_1, \cdots, t_{n_j}$ are full terms then $f_j(t_1, \cdots, t_{n_j})$
is a full term.
\end{enumerate}
\end{df}
\noindent
By $W_{\tau}^f(X)$ we denote the set of all full terms
of type $\tau$. It is easy to see that the set $W_{\tau}^f(X)$ of all full
terms is closed under composition.

\begin{lemma}\label{l2.2}
Let $s \in W_{\tau}^f(X_n)$ and let $t_1, \cdots, t_n \in W^f_{\tau}(X_m),
1 \leq n,m \in \natur$
be full terms. Then $S_m^n(s, t_1, \cdots, t_n)$ is also a full term.
\end{lemma}
{\bf Proof.} We give a proof by induction on the complexity ($Depth$) of
a term $s$.
If $s = f(x_{s(1)}, \cdots, x_{s(n)})$ where $s$
is a permutation on the set  $\{1, \cdots, n\}$ then
$S_m^n(s, t_1, \cdots, t_n) = f(t_{s(1)}, \cdots, t_{s(n)})$ is a full term by
Definition 2.1(ii).

Now, let $s = f(s_1, \cdots, s_r), s_1, \cdots, s_r \in W_{\tau}(X_n)$
and let  $s \in W_{\tau}^f(X_n),$ $ t_1, \cdots, t_n \in W_{\tau}^f(X_m)$.
We assume that $S_m^n(s_j, t_1, \cdots, t_n)$  are full terms for all
$1 \leq j \leq r$.  Then $$S_m^n(s, t_1, \cdots, t_n) =
f(S_m^n(s_1, t_1, \cdots, t_n), \cdots, S_m^n(s_r, t_1, \cdots, t_n))$$
and by Definition \ref{d2.1} (ii)  $S_m^n(s, t_1, \cdots, t_n)$ is a full term.
\hfill $\rule{2mm}{2mm}$

For full terms we have:
\begin{satz}
Let $s \in W_{\tau}^f(X_n)$ and assume that $t_1, \cdots,
t_n \in W_{\tau}(X_m), 1 \leq m,n \in \natur$ and that $\tau =(n, \cdots, n)$, i.e., all
operation symbols have the same arity. Then \\[-3mm]
\[Depth(S_m^n(s, t_1, \cdots, t_n)) = max\{Depth(t_1), \cdots, Depth(t_n)\}
+ Depth(s).\]
\end{satz}

{\bf Proof.} We give a proof by induction on $Depth(s)$.
If $Depth(s) = 1$ then $s = f(x_{s(1)}, \cdots, x_{s(n)})$ for an $n$-ary
operation symbol $f$ and a permutation $s:\{1, \cdots, n\} \to \{1, \cdots, n\}$.
There follows that $$Depth(S_m^n(s, t_1, \cdots, t_n)) =
Depth(f(t_{s(1)}, \cdots, t_{s(n)})) $$ $$= max\{Depth(t_{s(1)}), \cdots,
Depth(t_{s(n)})\} + 1 = max\{Depth(t_1), \cdots, Depth(t_n)\} + 1.$$

Assume that the formula is satisfied for $s_1, \cdots, s_r$ and
assume that $s = f(s_1, \cdots, s_r)$. Note that if  $s \in W_{\tau}^f(X_n)$
then also  $s_1, \cdots, s_r \in W_{\tau}^f(X_n)$, i.e. one can assume that
$$Depth(S_m^n(s_j, t_1, \cdots, t_n)) = max\{Depth(t_1), \cdots, Depth(t_n)\}
+ Depth(s_j)$$ for $1 \leq j \leq r$. Then we have
$$S_m^n(s, t_1, \cdots, t_n) = $$ $$S_m^n(f(s_1, \cdots, s_r), t_1, \cdots, t_n)
= f(S_m^n(s_1, t_1, \cdots, t_n), \cdots, S_m^n(s_r, t_1, \cdots, t_n))$$
and $$Depth(S_m^n(s, t_1, \cdots, t_n))$$ $$ =
max\{Depth(S_m^n(s_1, t_1, \cdots, t_n)), \cdots,
Depth(S_m^n(s_r, t_1, \cdots, t_n))\} +1 $$ $$= max\{max\{Depth(t_1), \cdots,
Depth(t_n)\} + Depth(s_1),$$ $$  \cdots, max\{Depth(t_1), \cdots,
Depth(t_n)\} + Depth(s_r)\}+1 $$ $$=
max\{Depth(t_1), \cdots, Depth(t_n)\} + max\{Depth(s_1), \cdots, Depth(s_r)\}
+1 $$ $$=  max\{Depth(t_1), \cdots, Depth(t_n)\} + Depth(s)$$since $s =
f(s_1, \cdots, s_r).$ \hfill $\rule{2mm}{2mm}$ \\

\section{The Depth of a Term with Respect to a Variable}

To derive a formula for the depth of the superposition of arbitrary terms we
define at first the depth of a term with respect to a variable.
Let $ t \in W_{\tau}(X_n)$ be an $n$-ary term and let $var(t)$ be the set of
all variables occurring in the term $t$.

\begin{df}~~
\begin{enumerate}
\item[\mbox{(i)}] If $t = x_k, 1 \leq k \leq n$, then $Depth_l(t):= 0$ for
all $1 \leq l \leq n$.
\item[\mbox{(ii)}] If $t = f_i(t_1, \cdots, t_{n_i})$ where $f_i$ is $n_i$-ary
and if we assume that $Depth_l(t_j), 1 \leq j \leq n_i, 1 \leq l \leq n$
are already defined then for all $l, 1 \leq l \leq n$ we define\\
\[
Depth_l(t):= \left \{ \begin{array}{ll}
0,& ~\mbox{if}~ x_l \not \in var(t)\\[1mm]
max\{Depth_l(t_j)|0 \leq j \leq n_i, x_l \in var(t_j)\}+1,&
 ~\mbox{otherwise.}
\end{array} \right. \]
\end{enumerate}
\end{df}

\begin{bsp}
Consider $I=\{1,2\}$ and  the type $\tau = (2,3)$.

For the term $t_1 = f_2(f_1(x_1,x_1), f_1(x_1,x_2), x_3)$ we have
$Depth_1(t_1) = 2, Depth_2(t_1) = 2, Depth_3(t_1) = 1$ and $Depth(t_1) = 2$.

For the term $t_2 = f_1(f_2(x_1,x_1,x_2),x_1)$ we have
$Depth_1(t_2) = 2, Depth_2(t_2) = 2, Depth_3(t_2) = 0$ and $Depth(t_2) = 2$.

For $t_3 = f_1(f_2(x_1,x_3,x_3), x_1)$ one has
$Depth_1(t_3) = 2, Depth_2(t_3) = 0, Depth_3(t_3) = 2$.

For $s = f_2(f_1(x_1,x_2), x_2, x_3)$ one obtains
$Depth_1(s) = 2, Depth_2(s) = 2$ and $Depth_3(s) = 1$.

Consider $S_3^3(s, t_1,t_2,t_3)$. Then it is easy to calculate that
$Depth(S_3^3(s, t_1,t_2,t_3)) = 4$. This is equal to $$
max\{Depth_1(s) + Depth(t_1),Depth_2(s) + Depth(t_2),
Depth_3(s) + Depth(t_3)\}.$$
\end{bsp}
More generally, we prove
\begin{satz}\label{t3.3}
Let $s \in W_{\tau}(X_n), t_1, \cdots, t_n \in W_{\tau}(X_m)$.
Then $$Depth(S_m^n(s, t_1, \cdots, t_n)) = max\{Depth_j(s) + Depth(t_j)|
1 \leq j \leq n, x_j \in var(s)\}.$$
\end{satz}
{\bf Proof.} We prove the formula by induction on $Depth(s)$.
If $Depth(s) = 0$  then there exists a natural number
$k \in \{1, \cdots, n\}$ such  that $s = x_k$ and then
$S_m^n(s, t_1, \cdots, t_n) = t_k$ and thus $Depth(S_m^n(s,t_1,\cdots, t_n))
= Depth(t_k)$ and $max\{Depth_j(s) + Depth(t_j)|1 \leq j \leq n, x_j \in var(s)\}
= Depth(x_k) + Depth(t_k) = 0 + Depth(t_k) = Depth(t_k)$.

Assume now that the formula is satisfied for $s_1, \cdots, s_r$ and assume that
$s=f(s_1, \cdots, s_r)$. \\Then $$S_m^n(s,t_1, \cdots, t_n) =
f(S_m^n(s_1, t_1, \cdots, t_n), \cdots,S_m^n(s_r, t_1, \cdots, t_n))$$ and
$$Depth(S_m^n(s,t_1, \cdots, t_n)) $$ $$=
max\{Depth(S_m^n(s_1, t_1, \cdots, t_n)), \cdots,
Depth(S_m^n(s_r, t_1, \cdots, t_n))\}+1$$ $$
=max\{max\{Depth_j(s_1) + Depth(t_j)|1 \leq j \leq n, x_j \in var(s_1)\}, \cdots,$$ $$
max\{Depth_j(s_r) + Depth(t_j)|1 \leq j \leq n, x_j \in var(s_r)\}\}+1$$ $$
= max\{max\{Depth_j(s_k)|1 \leq k \leq r, x_j \in var(s_k)\}+1 + Depth(t_j)|
1 \leq j \leq n, $$ $$x_j \in \bigcup\{var(s_k)|1 \leq k \leq r\}\}$$ $$
= max\{max\{Depth_j(s_k)|1 \leq k \leq r, x_j \in var(s_k)\}+1+Depth(t_j)|
1 \leq j \leq n,$$ $$ x_j \in var(s)\}
= max\{Depth_j(s) + Depth(t_j)| 1 \leq j \leq n,
x_j \in var(s)\}.$$  \hfill$\rule{2mm}{2mm}$\\

It is clear that the depth of a term $t$ is the maximum of all $Depth_j(t)$
for $x_j \in var(t)$, i.e.,

\begin{lemma}
If $t \in W_{\tau}(X_n), 1 \leq n \in \natur$ then
\[Depth(t) = max\{Depth_j(t)|1 \leq j \leq n, x_j \in var(t)\}.\]
\end{lemma}
{\bf Proof.} We will give a proof by induction on $Depth(t)$.
If $Depth(t) = 0$ then $t = x \in X$ is a variable and $Depth_j(t) = 0$
for $1 \leq j \leq n$. Hence $max\{Depth_j(t)|1 \leq j \leq n, x_j
\in var(t)\} = 0 = Depth(t)$.

Assume the lemma is satisfied for $s_1, \cdots, s_r$ and that
$t = f(s_1, \cdots,s_r)$.
Then $
Depth(t)\\ = max\{Depth(s_k)|k \in \{1, \cdots, r\}\}+1 \\
= max\{max\{Depth_j(s_k)|1 \leq j \leq n, x_j \in var(s_k)\}|1
\leq k \leq r\}+1\\
= max\{max\{Depth_j(s_k)|1 \leq k \leq r, x_j \in var(s_k)\}+1 |1
\leq j \leq n, x_j \in var(t)\}\\
= max\{Depth_j(t)|1 \leq j \leq n, x_j \in var(t)\}$. \hfill $\rule{2mm}{2mm}$\\

\section{Full Hypersubstitutions}

In section 1 we have already introduced  the concept of a hypersubstitution.
Let $Hyp(\tau)$ be the set of all hypersubstitutions of type $\tau$. On
$Hyp(\tau)$ by $(\sigma_1 \circ_h \sigma_2)(f_i) :=
\hat \sigma_1[\sigma_2(f_i)]$ for all operation symbols $f_i$ a binary
operation
can be defined. Then $Hyp(\tau)$ together with the identity hypersubstitution
$\sigma_{id}$ defined by $\sigma_{id}(f_i) := f_i(x_1, \cdots, x_{n_i})$
forms a monoid. If $V$ is a variety of algebras of type $\tau$ then $V$
is called solid if for every identity $s \approx t$ in $V$ and every
$\sigma \in  Hyp(\tau)$ the equations $\hat \sigma[s] \approx \hat \sigma[t]$
are satisfied as identities in $V$. All solid varieties of type $\tau$
form a complete sublattice of the lattice of all varieties of type $\tau$.
If $M$ is a submonoid of $Hyp(\tau)$ one can define $M$-solid varieties of
type $\tau$. The class of all $M$-solid varieties of type $\tau$ forms also
a complete lattice and the collection of all solid varieties is a complete
sublattice of the lattice of all $M$-solid varieties of type $\tau$. More generally,
if $M_1 \subseteq M_2$ for two submonoids $M_1, M_2$ of $Hyp(\tau)$
then the collection of all $M_2$-solid varieties of type $\tau$ forms a complete
sublattice of the lattice of all $M_1$-solid varieties of type $\tau$.

For more background on hyperidentities, hypersubstitutions and solid
varieties see \cite{Den1,Den2}.

Now we consider a special class of hypersubstitutions of type $\tau$.
\begin{df}
A hypersubstitution is called full if $\sigma(f_i) \in W_{\tau}^f(X_{n_i})$
for all $i \in I$. By $Hyp^f(\tau)$ we denote the set of all full
hypersubstitutions of type $\tau$.
\end{df}
\begin{lemma}
The set $Hyp^f(\tau)$ forms a submonoid of the monoid $Hyp(\tau)$ of all
hypersubstitutions of type $\tau$.
\end{lemma}
{\bf Proof.} Since the terms $\sigma_{id}(f_i) = f_i(x_1, \cdots, x_{n_i})$
are  full terms for every $i \in I$ the identity hypersubstitution is full.
Assume that $\sigma_1, \sigma_2 \in Hyp^f(\tau)$. We want to prove that
$(\sigma_1 \circ_h \sigma_2)(f_i)$ are full terms for every $i \in I$.
Since $\sigma_2(f_i)$ is a full term, by definition of full terms there exists
an operation symbol $f_j$ such that $\sigma_2(f_i) = f_j(x_{s(1)}, \cdots,
x_{s(x_{n_j})})$ for a permutation $s: \{x_1, \cdots, x_{n_j}\} \to
\{x_1, \cdots, x_{n_j}\}$  or full terms $t_1, \cdots t_{n_j}$ with
$\sigma_2(f_i) = f_j(t_1, \cdots, t_{n_j})$.
In the first case we have $$(\sigma_1 \circ_h \sigma_2)(f_i)=
\hat \sigma_1[\sigma_2(f_i)]=\sigma_1(f_j)(x_{s(1)}, \cdots, x_{s(n_j)}).$$

Since $\sigma_1$ is a full
hypersubstitution, the term $\sigma_1(f_{j})$ is a full term and then for
every permutation $s:\{1, \cdots, n_j\} \to \{1, \cdots, n_j\}$ the term
$\sigma_1(f_j)(x_{s(1)}, \cdots, x_{s(n_j)})$ is also full.
In the second case one obtains
$$(\sigma_1 \circ_h \sigma_2)(f_i)=
\hat \sigma_1[\sigma_2(f_i)]=\hat \sigma_1[f_j(t_1, \cdots, t_{n_j})]
= \sigma_1(f_j)(\hat\sigma_1[t_1], \cdots, \hat \sigma_1[t_{n_j}]).$$

Here $\sigma_1(f_j)$ is a full term. We show that $\hat \sigma[t_k]$
are also full terms for all $k \in \{1, \cdots, n_j\}$.
In fact, if $t_k = f_{\mu}(x_{s(1)}, \cdots, x_{s(n_{\mu})})$
then $\hat \sigma_1[t_k]
= \sigma_1(f_{\mu})(x_{s(1)}, \cdots, x_{s(n_{\mu})})$
where $s:\{1, \cdots, n_{\mu}\}
\to \{1, \cdots, n_{\mu}\}$ is a permutation. Since $\sigma_1$ is a full
hypersubstitution, the term $\sigma_1(f_{\mu})$ is a full term and then for
every permutation $s:\{1, \cdots, n_{\mu}\} \to \{1, \cdots, n_{\mu}\}$ the term
$\sigma_1(f_{\mu})(x_{s(1)}, \cdots, x_{s(n_{\mu})})$ is also full. If $t_k =
f_{\mu}(t_{11}, \cdots, t_{1n_{\mu}})$ and assume that
$\hat \sigma_1[t_{1j}],
1 \leq j \leq n_{\mu}$ are full then $\hat \sigma_1[t_k]
= \sigma_1(f_{\mu})(\hat
\sigma[t_{11}], \cdots, \hat \sigma[t_{1n_{\mu}}])$ is also full by Definition \ref{d2.1} (ii).
But then by Lemma \ref{l2.2} $\sigma_1(f_j)(\hat\sigma[t_1], \cdots
\hat \sigma[t_{n_j}])$ is also a full term and
$\sigma_1 \circ_h \sigma_2 \in Hyp^f(\tau)$. \hfill $\rule{2mm}{2mm}$ \\

For a given variety $V$ of type $\tau$ one can determine all subvarieties
of $V$ which are $Hyp^f(V)$-solid. Consider as an example the variety
of all bands (idempotent semigroups).
A hypersubstitution $\sigma$ of type $\tau$ is called a regular
hypersubstitution of type $\tau$ if $var(\sigma(f_i)) =
\{x_1, \cdots, x_{n_i}\}$ for all $i \in I$. The collection $Reg(\tau)$
of all regular hypersubstitutions of type $\tau$ forms also a
submonoid $Reg(\tau)$ of $Hyp(\tau)$ and a variety $V$ of type
$\tau$ is called regular-solid if it is $Reg(\tau)$-solid (\cite{Den2}).\\

Then we get
\begin{prop}
A variety $V$ of bands is $Hyp^f(2)$-solid iff $V$ is regular-solid.
\end{prop}\noindent
{\bf Proof.} By definition every full hypersubstitution is regular, i.e.,
$Hyp^f(2) \subseteq Reg(2)$.
Therefore, if $V$ is regular-solid then it is also $Hyp^f(2)$-solid.
It remains to show that a $Hyp^f(2)$-solid variety of bands is regular-solid.
Let $\sigma_t \in Reg(2)$ where $t \in W_2(X_2)$. (Here $\sigma_t$ is the
hypersubstitution which maps the binary operation symbol $f$ to the
binary term $t$.) In $t$ both variables $x_1$ and $x_2$  occur
and there are exactly the following cases:

a)  $t$ starts with $x_1$ and ends with $x_2$,

b)  $t$ starts with $x_2$ and ends with $x_1$,

c) $t$ starts and ends with $x_1$,

d) $t$ starts and ends with $x_2$.

In the first case we have $t \approx x_1x_2$, in the second case there holds
$t \approx x_2x_1$, in the third case one obtains
$t \approx x_1x_2x_2x_1$ and in the
fourth case $t \approx x_2x_1x_1x_2$.
All four hypersubstitutions are full and since $V$ is $Hyp^f(2)$-solid
it is also regular-solid. ~~~~~~~\hfill $\rule{2mm}{2mm}$ \\

\noindent
If the depth of $\sigma(f_i)$ for a hypersubstitution $\sigma$ for all
$i \in I$ is known and if  $Depth(t)$ for $t \in W_{\tau}(X)$ is known
we want to know what $Depth(\hat \sigma[t])$ is.
It is quite natural to define the depth of a hypersubstitution
$\sigma$ in the  following way:

\begin{df}
Let $\sigma$ be a hypersubstitution of type $\tau$ then
\[Depth(\sigma) := max\{Depth(\sigma(f_i))|i \in I \}.\]
\end{df}

Clearly, for the type $\tau = (n)$ we have
$Depth(\sigma) = Depth(\sigma(f))$.
Then we obtain the following result:
\begin{cor}\label{c4.5}
Let $t \in W^f_{\tau}(X_n), \sigma \in Hyp^f(\tau)$ and $\tau = (n)$ for a
natural number $n \geq 1$. Then
\[Depth(\hat \sigma[t]) = Depth(\sigma(f))  Depth(t).\]
\end{cor}
{\bf Proof.} We give a proof by induction on $Depth(t)$. If $Depth(t) = 1$
then $t = f(x_{s(1)}, \cdots, x_{s(n)})$ for a permutation $s:\{1, \cdots, n\}
\to \{1, \cdots, n\}$ and $\hat \sigma[t] = \sigma(f)(x_{s(1)}, \cdots,
x_{s(n)}) = S_n^n(\sigma(f), x_{s(1)}, \cdots, x_{s(n)})$. Therefore\\
$\begin{array}{ll} Depth(\hat \sigma[t])&= Depth(S_n^n(\sigma(f), x_{s(1)}, \cdots, x_{s(n)}))\\
&= max\{Depth(x_{s(1)}), \cdots, Depth(x_{s(n)})\} + Depth(\sigma(f))\\ &= 0 +
Depth(\sigma(f)) \\ &= Depth(\sigma(f))
= Depth(\sigma(f))  Depth(t).\end{array}$

Assume that $Depth(\hat\sigma[s_k])) = Depth(\sigma(f))  Depth(s_k)$
for $1 \leq k \leq n$ and $t = f(s_1, \cdots, s_n)$.
Then\\
 $\begin{array}{ll} Depth(\hat \sigma[t]) &= Depth(S_n^n(\sigma(f), \hat\sigma[s_1], \cdots,
\hat \sigma[s_n])) \\ &
= Depth((\sigma(f)) + max\{Depth(\hat\sigma[s_k])|1 \leq k \leq n\}\\
&= Depth((\sigma(f)) + max\{Depth(\sigma(f))  Depth(s_k)|1 \leq k \leq n\}\\
&= Depth((\sigma(f)) + Depth(\sigma(f)) max\{Depth(s_k)|1 \leq k \leq n\}\\
&= Depth((\sigma(f)) (1+max\{Depth(s_k)|1 \leq k \leq n\})\\
&= Depth(\sigma(f))  Depth(t).\end{array}$

 ~~~~~~~~~~~~ \hfill $\rule{2mm}{2mm}$ \\

As a consequence of Corollary \ref{c4.5} we have:

\begin{cor}
The function $Depth^f: Hyp^f(\tau) \to  \natur^+$ with
$\sigma \mapsto Depth(\sigma)$ defines a homomorphism
from the monoid $(Hyp^f(\tau); \circ_h, \sigma_{id})$ onto the monoid
$(\natur^+; \cdot, 1)$.
\end{cor}
{\bf Proof.} The mapping $Depth^f$ is well-defined and surjective since to every
natural number $n \geq 1$ there is a full term $t$ with $Depth(t) =n$.
Further we have\\
 $\begin{array}{ll} Depth^f(\sigma_{id}) &= Depth(\sigma_{id}(f)) \\ &=
Depth(f(x_1, \cdots, x_n)) = 1\end{array}$ \\
and \\
$\begin{array}{ll} Depth^f(\sigma_1 \circ_h \sigma_2) &=
Depth((\sigma_1 \circ_h \sigma_2)(f))\\ & =
Depth(\hat\sigma_1[\sigma_2(f)]) \\ &=
Depth(\sigma_1(f))  Depth(\sigma_2(f)) \\ &= Depth^f(\sigma_1)
Depth^f(\sigma_2). \end{array}$

~~~~ \hfill $\rule{2mm}{2mm}$

\section{The Depth of an Arbitrary Hypersubstitution}

To derive a formula for $Depth(\hat \sigma[t])$ we introduce the following notation:

The yield $yd(t)$ of the term $t$ is defined inductively as follows:
\begin{enumerate}
\item[\mbox{(i)}] $yd(x) = x$ for all $x \in X$.
\item[\mbox{(ii)}] If $t = f_i(t_1, \cdots, t_{n_i})$ then $yd(t) = yd(t_1) \cdots
yd(t_{n_i})$.
\end{enumerate}
That means, the yield of a term $t$ is the semigroup word obtained from $t$ by
cancellation of all operation symbols, brackets, and commas.

The length $\ell(t)$ of the term $t$ is the number of variables occurring in $t$.

Assume that $t \in W_{\tau}(X)$ has the length $\ell(t) = n$ and that
$i \in \{1, \cdots, n\}.$
Let $yd(t) = u_1 \cdots u_{\ell(t)}$, where $u_1 \cdots u_{\ell(t)}$ denote certain variables. Then we define terms $A_k(i,t)$ for
$0 \leq k \in \natur$ and for $1 \leq i \leq \ell(t)$ in the following
inductive way:
\begin{enumerate}
\item[\mbox{(i)}] $A_0(i,t):= u_i$,
\item[\mbox{(ii)}]  $A_{\ell(t)}':= t$,
\item[] $A_{k}':= t_r$, if $A_{k+1}':= f_j(t_1,\ldots t_{n_j})$  and $u_i$ is contained as a subterm in $t_r$ $(1\leq r\leq n_j)$ (where $r$ is uniquely determined),
\item[] $A_{k}':= u_i$, if $A_{k+1}':= u_i$  for $0\leq k <\ell(t)$,
\item[\mbox{(iii)}] Let $r$ be the greatest integer with $A_{r}'= u_i$. Then we define
\item[] $A_{i}:=A_{i+r}'$ for $0\leq i <\ell(t)-r$ and
\item[] $A_{i}:=t$ for $\ell(t)-r<k\in\natur$.
\end{enumerate}

By $\beta(i,t)$ we denote the least natural number $k$ with $A_k(i,t) = t$.
Let $\sigma \in Hyp(\tau)$. Then we define:
\begin{enumerate}
\item[\mbox{(i)}] $B_0(\sigma, i, t) := 0$,
\item[\mbox{(ii)}]  $B_k(\sigma, i, t) := Depth_a(\sigma(f_j))$ if $1 \leq
k \leq \beta(i,t)$ and $a \in \{1,2, \cdots, n_j\}$ is determined by the
property that $A_{k-1}(i,t)$ is at the $a$-th place in the term
$A_k(i,t) = f_j(t_1, \cdots, t_{n_j})$.
\end{enumerate}
Finally we define
$B(\sigma, i,t):= \sum \limits_{k=0}^{\beta(i,t)} B_k(\sigma, i,t)$
and\\
$B(\sigma, t) := max(\{B(\sigma,i,t)|1 \leq i \leq \ell(t), B_{\beta(\sigma,i,t)}(i,t) \not
= 0\} \cup \{0\})$

Then we have
\begin{satz}
Let $t \in W_{\tau}(X)$ and $\sigma \in Hyp(\tau)$. Then $Depth(\hat\sigma[t]) =
B(\sigma, t)$.
\end{satz}

{\bf Proof.} We will give a proof by induction on $Depth(t)$.
If $Depth(t) = 0$ then $t \in X, \ell(t) = 1, A_0(1,t) = t, $ and
$\beta(1,t) = 0$. Then follows  $$B(\sigma,t) =
max(\{B(\sigma, i, t)| 1 \leq i \leq \ell(t), B_{\beta(\sigma,i,t)}(i,t) \not = 0\}
\cup \{0\}) $$ $$= 0 = Depth(t) = Depth(\hat \sigma[t]).$$

Assume that $Depth(\hat \sigma[t]) = B(\sigma, t)$ if $Depth(t) < p$
for a natural number $p \geq 1$. If $Depth(t) = p$ then there are terms
$t_1, \cdots, t_{n_i}$ and an $n_i$ary operation symbol $f_i$
with $t=f(t_1, \cdots, t_{n_i})$ and $Depth(t_k)<p$ for $k =1, \cdots, n_i$.
Then we have
\\
$\begin{array}{ll} Depth(\hat\sigma[t]) &=
Depth(S_n^{n_i}(\sigma(f_i), \hat\sigma[t_1], \cdots, \hat\sigma[t_{n_i}])) \\
&= max\{Depth_j(\sigma(f_i)) + Depth(\hat \sigma[t_j]) |1 \leq j \leq n_i,
x_j \in var(\sigma(f_i)) \}  \\ & \mbox{(by Theorem \ref{t3.3})} \\
&= max\{Depth_j(\sigma(f_i)) + B(\sigma, t_j) |1 \leq j \leq n_i,
x_j \in var(\sigma(f_i))\}\\ & \mbox{(by our hypothesis)} \hfill (*) \\
&= max\{Depth_j(\sigma(f_i)) + max(\{B(\sigma, i, t_j)|
1 \leq i \leq \ell(t_j), \\ &B_{\beta(i,t_j)}(\sigma,i,t_j) \neq 0 \} \cup \{0\})\}|
1 \leq j \leq n_i, x_j \in var(\sigma(f_i)))\}.\end{array}$

If $n(j), 1 \leq j \leq n_i$ denotes the length of $t_j$, then we define:
$$q:= \sum \limits_{j=1}^{n_i} n(j), \ \
\delta_j:= \sum\limits_{a=1}^{j}n(a)\ \ \mbox{and}\ \ \beta_j = q - \sum
\limits_{j=1}^{n_i}n(a).$$
It is easy to check that $\beta(i,t) = \beta(i-\delta_j,t_j) +1$ for
$1 \leq j \leq n_i, \delta_j < i \leq \beta_j$ and $B_a(\sigma,i,t) =
B_a(\sigma, i - \delta_j, t_j)$ for $1 \leq j \leq n_i,
\delta_j < i \leq \beta_j$ and $0 \leq a \leq \beta(i-\delta_j,t_j)$.
Since $$A_{\beta(i,t)}(i,t) = t = f_i(t_1, \cdots, t_{n_i}),\ \
 A_{\beta(i- \delta_j, s_j)}(i,t) = t_j$$ for $\delta_j < i \leq \beta_j$ and
$1 \leq j \leq n_i$ we have $B_{\beta(i,t)}(\sigma,i,t) = Depth_j(\sigma(f_i))$ for
$\delta_j < i \leq \beta_j$ and $1 \leq j \leq n_i$.
Hence for $1 \leq j \leq n_i$ and $1 \leq i \leq n(j)$ there holds
$$B(\sigma, i, t_j) = \sum \limits_{a=0}^{\beta(i,t_j)}B_a(\sigma,i,t_j)
= \sum\limits_{a=0}^{\beta(i+\delta_j,t)-1}B_a(\sigma, i+\delta_j,t).$$

Now we continue with (*) and obtain by substitution
$$Depth(\hat\sigma[t])= $$ $$= max(\{B_{\beta(i+\delta_j,t)}(\sigma, i+\delta_j,t) +
max(\{\sum\limits_{a=0}^{\beta(i+\delta_j,t) -1}B_a(\sigma, i+\delta_j,t)|
1 \leq i \leq n(j), $$ $$B_{\beta(i + \delta_j,t)-1}(\sigma,i + \delta_j, t) \neq 0\}
\cup \{0\})|1 \leq j \leq n_i, x_j \in var(\sigma(f_i))\} $$ $$
= max(\{B_{\beta(i+\delta_j,t)}(\sigma, i+\delta_j,t) + \sum\limits_{a=0}
^{\beta(i+ \delta_j,t)-1}B_a(\sigma, i+ \delta_j,t)$$ $$
|1 \leq i \leq n(j), 1 \leq j \leq n_i, x_j \in var(\sigma(f_i)),
B_{\beta(i+\delta_j,t)-1}(\sigma,i + \delta_j,t) \not = 0 \} \cup \{0\})$$
(since $B_{\beta(i+\delta_j,t)}(\sigma, i+\delta_j,t) =
const_j$ for $1 \leq i \leq n(j)$ and $1 \leq j \leq n_i$ for a natural number
$const_j$)
$$= max(\{\sum\limits_{a=0}^{\beta(i+\delta_j,t)}B_a(\sigma, i + \delta_j,t)
|1 \leq i \leq n(j), 1 \leq j \leq n_i, $$ $$B_{\beta(i + \delta_j,t)-1}
(\sigma,i + \delta_j,t) \not = 0\} \cup \{0\})$$ $$
= max(\{B(\sigma, i+ \delta_j,t)|1\leq i \leq n(j), 1 \leq j \leq n_i,
B_{\beta(i + \delta_j,t)-1}
(\sigma,i + \delta_j,t) \not = 0\} \cup \{0\})$$ $$
= max(\{B(\sigma,i,t)|1 \leq i \leq \ell(t), B_{\beta(i,t)}(i,t) \not = 0\}
\cup\{0\})
= B(\sigma,t).$$\hfill $\rule{2mm}{2mm}$\\

\vspace*{3mm}

K. Denecke\\
Universit\"at Potsdam\\
Fachbereich Mathematik\\
Postfach 601553\\
D-14415 Potsdam\\
email: {\ttfamily kdenecke@rz.uni-potsdam.de}

J. Koppitz\\
Universit\"at Potsdam\\
Fachbereich Mathematik\\
Postfach 601553\\
D-14415 Potsdam\\
email: {\ttfamily koppitz@rz.uni-potsdam.de}

Sl. Shtrakov\\
South-West-University Blagoevgrad\\
Faculty of Mathematics and Natural Sciences\\
2700 Blagoevgrad\\
Bulgaria\\
e-mail:{\ttfamily shtrakov@swu.bg}

\end{document}